\documentclass{article}
\usepackage[english]{babel}
\usepackage{amsmath,amssymb,textcomp,bbm,epsfig}

\newcommand{\assign}{:=}
\newcommand{\nosymbol}{}
\newcommand{\tmop}[1]{\ensuremath{\operatorname{#1}}}
\newcommand{\tmstrong}[1]{\textbf{#1}}
\newcommand{\tmtextbf}[1]{{\bfseries{#1}}}
\begin{document}

\title{Decomposing manifolds into Cartesian products\footnote{2010 Mathematics Subject
Classification Primary 57R80. Secondary 57M50. Key words and phrases: Seifert
manifolds, Whitehead torsion, s-cobordisms, surgery groups. The first author
acknowledges the support of the Simons Foundation Grant No. 281810.}}

\author{Slawomir Kwasik \\
Department of Mathematics, Tulane University, New Orleans, LA 70118, USA\\
kwasik@tulane.edu\\
\\
AND Reinhard Schultz\\
Department of Mathematics, University of California, Riverside, CA 92521, USA\\
schultz@math.ucr.edu}

\date{November 22, 2016}

\maketitle

\begin{abstract}
  The decomposability of a Cartesian product of two nondecomposable manifolds into
  products of lower dimensional manifolds is studied.  For $3$-manifolds we obtain
  an analog of a result due to Borsuk for surfaces, and in higher dimensions we show
  that similar analogs do not exist unless one imposes further restrictions such as
  simple connectivity.
\end{abstract}

\section{Introduction}

There are plenty of examples which show the nonuniqueness of splitting a space
(manifold) into Cartesian products. For example, there is the well known Bing space,
a generalized 3-manifold $X$ ({\it cf.} [Bi]) such that $X \neq
\mathbb{R}^3$ and $X \times \mathbb{R} = \mathbb{R}^4 = \mathbb{R}^3
\times \mathbb{R}$, or the open 3-manifold $\mathcal{W}$ of J. H. C.
Whitehead ({\it cf.} [He]), where $\mathcal{W} \neq \mathbb{R}^3$ and
again $\mathcal{W} \times
\mathbb{R} = \mathbb{R}^4 = \mathbb{R}^3 \times \mathbb{R}$ (here =
stands for homeomorphic). More dramatic examples are pairs of two
homotopy inequivalent 3-dimensional Seifert manifolds $\mathcal{M}$ and
$\mathcal{N}$ such that $\mathcal{M} \times S^1 = \mathcal{N} \times
S^1$ ({\it cf.} [CR],[KR2]).

However,
on the positive side there is an old result of K.Borsuk ({\it cf.} [Bo])
that a closed, n-dimensional manifold has at most one decomposition into the
Cartesian product of indecomposable factors of dimension $\leqslant 2$.

Now suppose we that have two closed, oriented $n$-dimensional manifolds
$\mathcal{M}^n$ and $\mathcal{N}^n$ which can not be split into products of
closed, oriented manifolds ($\neq \{\tmop{pt}\}$) of lower dimension.  Here is
one natural question: {\sl Can $\mathcal{M}^n \times \mathcal{N}^n$ be
decomposed into products of manifolds of dimension $\leqslant n - 1$?}\\

More generally, we shall consider the following situation(we recall hat we are working in the TOP category, i.e., topological manifolds and homeomorphisms):

Let $\mathcal{M}^n$, $\mathcal{N}^k$ be closed, oriented, indecomposable
(into
nontrivial Cartesian products) manifolds of dimension $n$ and $k$
respectively, $k \leqslant n$. One says that {\tmstrong{the manifold}}
$\mathcal{N}^k$ {\tmstrong{stably decomposes}} $\mathcal{M}^n$ if
$\mathcal{M}^n \times \mathcal{N}^k$ can be written as a Cartesian product of
manifolds of dimension $\leqslant n - 1$ (i.e. $\mathcal{M}^n \times
\mathcal{N}^k = \underset{i}{\Pi}~Y_i^{n_i}$ such that each $Y_i^{n_i}$ is a
closed manifold of dimension $n_i$, where $1 \leqslant n_i \leqslant n - 1$ and
$\underset{i}{\Sigma}\;n_i = n + k$). If for a given manifold $\mathcal{M}^n$,
there is no such $\mathcal{N}^k$, then $\mathcal{M}^n$ is called
{\tmstrong{stably nondecomposable}}.

If $n = 1$ or $n = 2$, then Borsuk's result shows that every $\mathcal{M}^n$ is
stably nondecomposable. It turns out that this is also true for $n = 3$:\\

\underline{{\tmstrong{Theorem A}}}: {\tmstrong{Let $\mathcal{M}^3$ be an
oriented, closed, nondecomposable 3-manifold. Then $\mathcal{M}^3$ is stably
nondecomposable.}}\\

On the other hand, for $n = 4$ we have the following:\\

{\tmstrong{\underline{Theorem B}: There exists an oriented, closed,
nondecomposable 4-manifold $\mathcal{M}^4$ such that $\mathcal{M}^4 \times S^k = S^1 \times S^k \times \tmop{RP}^3~(k =
2, 3, 4)$. Moreover there are infinitely many non-decomposable 4-manifolds $\mathcal{M}^4_i (i=1,2,...,n,...)$ with $\mathcal{M}^4_i \neq \mathcal{M}^4_j,i \neq j$ and $\mathcal{M}^4 \times S^k = S^1 \times S^k \times RP^3 (k=2,3,4)$.}}

\

The manifold $\mathcal{M}^4$ in the above theorem is not simply connected. It
turns out that this is an essential condition in our proof. Indeed, for simply
connected 4-manifolds we
have the following addendum to Theorem B.\\

{\tmstrong{\underline{Theorem B'}: Let $\mathcal{M}^4$ be a closed, simply
connected nondecomposable manifold. Then $S^k~(k = 2, 3, 4)$ cannot stably
decompose $\mathcal{M}^4$.}}

\section{Proofs}

This section contains proofs of our results. The methods and techniques
employed in these proofs form a curious combination of high-dimensional
surgery theory and low-dimensional topology.\\

\underline{{\tmstrong{Proof of Theorem A}}}: We first consider the case of
$\mathcal{M}^3$ and $\mathcal{N}^3$. Let $\mathcal{M}^3$ and $\mathcal{N}^3$
be oriented, closed, nondecomposable 3-manifolds.

Suppose $\mathcal{M}^3 \times \mathcal{N}^3$ is decomposable, so we can write
$\mathcal{M}^3 \times \mathcal{N}^3 = S_1 \times S_2 \times S_3$, where
$\tmop{dim}\,S_i = 2 (i = 1, 2, 3)$.

Our first observation is that without loss of generality we can assume $\pi_1
(\mathcal{M}^3)$, $\pi_1 (\mathcal{N}^3)$ are infinite. Our second observation
is that because the Euler characteristic $\chi (\mathcal{M}^3) = \chi
(\mathcal{N}^3) = 0$, then at least one of the $S_i (i = 1, 2, 3)$ must be a
torus $T^2 = S^1 \times S^1$.

Suppose that exactly one of the $S_i$ (say $S^{\nosymbol}_3$) is a torus. Then
$\mathcal{M}^3 \times \mathcal{N}^3 = S_1 \times S_2 \times S^1 \times S^1$.
It follows that the center $C (\pi_1 (\mathcal{M}^3 \times \mathcal{N}^3))$ of
$\pi_1  (\mathcal{M}^3 \times \mathcal{N}^3)$ is given by $C (\pi_1
(\mathcal{M}^3 \times \mathcal{N}^3)) \cong C (\pi_1 (\mathcal{M}^3)) \oplus C
(\pi_1 (\mathcal{N}^3)) \cong \mathbb{Z} \oplus \mathbb{Z}$.

Now suppose that the center of $\pi_1 (\mathcal{M}^3)$ or the center of \
$\pi_1 (\mathcal{N}^3)$ is trivial, say $C (\pi_1 (\mathcal{M}^3)) \cong 0$.
Consequently, $C (\pi_1 (\mathcal{N}^3)) \cong \mathbb{Z} \oplus
\mathbb{Z}$. This implies that $\mathcal{N}^3 = S^1 \times S^1 \times
S^1$ ({\it i.e.\/}, Theorem 12.10, p.131 in [He]), and we are done.

Assume now that $C (\pi_1 (\mathcal{M}^3)) \cong C (\pi_1 (\mathcal{N}^3))
\cong \mathbb{Z}$. This implies that both $\mathcal{M}^3$ and $\mathcal{N}^3$
are Seifert manifolds ({\it cf.} [CJ], [G]). Let $h : \mathcal{M}^3 \times
\mathcal{N}^3 \overset{\approx}{\rightarrow} S_1 \times S_2 \times S^1 \times
S^1$ be a homeomorphism. The induced homomorphism
\[ h_{\ast} : \pi_1  (\mathcal{M}^3 \times \mathcal{N}^3) \rightarrow \pi_1
   (S_1 \times S_2 \times S^1 \times S^1) \]
is an isomorphism and
\[ h_{\ast |} : C (\pi_1 (\mathcal{M}^3 \times \mathcal{N}^3)) \rightarrow C
   (\pi_1 (S_1 \times S_2 \times S^1 \times S^1)) \]
is an isomorphism as well, i.e.,
\[ h_{\ast |} : \mathbb{Z} \oplus \mathbb{Z}
   \overset{\cong}{\longrightarrow} \mathbb{Z} \oplus \mathbb{Z} \]
To go further we resort to the following simple torus trick:
Since every automorphism $\mathbb{Z} \oplus \mathbb{Z}
\overset{\cong}{\longrightarrow} \mathbb{Z} \oplus \mathbb{Z}$ can be
realized by a homeomorphism $S^1 \times S^1 \rightarrow S^1 \times S^1$, that
is, $\pi_0 (\tmop{Homeo} (T^2)) \cong \tmop{GL} (2, \mathbb{Z})$ ({\it e.g.\/}, see
Theorem 4, p.26 in [R]). Then by composing
\[ h : \mathcal{M}^3 \times \mathcal{N}^3 \longrightarrow S_1 \times S_2
   \times S^1 \times S^1 \]
with $h' = \tmop{id}_{S_1 \times S_2} \times f : S_1 \times S_2 \times S^1
\times S^1 \longrightarrow S_1 \times S_2 \times S^1 \times S^1$ for some
homeomorphism $f : S^1 \times S^1 \rightarrow S^1 \times S^1$, we can assume
that there is a homeomorphism
\[ h : \mathcal{M}^3 \times \mathcal{N}^3 \longrightarrow S_1 \times S_2
   \times S^1 \times S^1 \]
with the isomorphism $h_{\ast |} : \mathbb{Z} \oplus \mathbb{Z}
\overset{\cong}{\longrightarrow} \mathbb{Z} \oplus \mathbb{Z}$ given by
\[ h_{\ast |} : h'_{\ast} \oplus h'_{\ast}~. \]
This implies that there is an induced homeomorphism
\[ \tilde{h} : \widetilde{\mathcal{M}^3} \times \widetilde{\mathcal{N}^3}
   \longrightarrow S_1 \times S_2 \times \mathbb{R} \times \mathbb{R} \]
where $\widetilde{\mathcal{M}^3}$, $\widetilde{\mathcal{N}^3}$ are infinite
cyclic coverings determined by the corresponding centers.

The fundamental groups $\pi_1 (\widetilde{\mathcal{M}^3}) \cong \pi_1
(\mathcal{M}^3) / C (\pi_1 (\mathcal{M}^3))$ and $\pi_1
(\widetilde{\mathcal{N}^3}) \cong \pi_1 (\mathcal{N}^3) / C (\pi_1
(\mathcal{N}^3))$ are Fuchsian groups (since both $\mathcal{M}^3$,
$\mathcal{N}^3$ are Seifert manifolds) and there is an isomorphism
\[ \pi_1 (\widetilde{\mathcal{M}^3}) \times \pi_1 (\widetilde{\mathcal{N}^3})
   \cong \pi_1 (S_1) \times \pi_1 (S_2) \]
This implies (by Proposition II.37, p.19, in [JS]) that $\pi_1 (\widetilde{M^3})$ and
$\pi_1 (\widetilde{\mathcal{N}^3})$ are isomorphic to fundamental groups of
closed surfaces. It is not difficult to see (for example using the cohomology
ring structure of closed surfaces) that corresponding groups must be
isomorphic, say
\[ \pi_1 (\widetilde{\mathcal{M}^3}) \cong \pi_1 (S_1)\quad \tmop{and}\quad \pi_1
   (\widetilde{\mathcal{N}^3}) \cong \pi_1 (S_2) \]
In particular
\[ \widetilde{\mathcal{M}^3} \times \mathcal{N}^3 = S_1 \times S_2 \times \mathbb{R}
   \times S^1 \]
and since $\pi_1 (\widetilde{\mathcal{M}^3}) \cong \pi_1 (S_1)$ then we have a
homotopy equivalence
\[ S_1 \times \mathcal{N}^3 \simeq S_1 \times S_2 \times S^1 \]
and by symmetry
\[ \mathcal{M}^3 \times \widetilde{\mathcal{N}^3} = S_1 \times S_2 \times S^1
   \times \mathbb{R} \]
which gives
\[ \mathcal{M}^3 \times S_2 \simeq S_1 \times S_2 \times S^1 \]
If one of the $S_i, i = 1, 2$ is $S^2$, say $S_1 = S^2$, then $S^2 \times
\mathcal{N}^3 \simeq S^2 \times S_2 \times S^1$ and
\[ \pi_1 (\mathcal{N}^3) \cong \pi_1  (S_2 \times S^1) \cong \pi_1 (S_2)
   \oplus \pi_1 (S^1) \]
Consequently $\mathcal{N}^3 = S_2 \times S^1$ ({\it cf.} [He] p.114). In particular
we can assume $S_1 \neq S^2, S_2 \neq S^2$.

Going back to the homotopy equivalence
\[ S_1 \times \mathcal{N}^3 \simeq S_1 \times S_2 \times S^1 \]
we observe that $S_1 \times S_2 \times S^1$ admits a metric (the standard
product metric) of non-positive curvature. The results of Farrell-Jones
({\it cf.} [FJ]) imply $S_1 \times \mathcal{N}^3 = S_1 \times (S_2 \times S^1)$. We
claim that this is impossible given the indecomposability of $\mathcal{N}^3$.

To see this we need the following slight adjustment in the conclusion of
Theorem 1 in [KR].\\

{\underline{{\tmstrong{Claim}}}}: Let $X, Y$ be closed oriented surfaces of
genus at least 2, and $\mathcal{N}^3$ be a Seifert manifold. If $\mathcal{N}^3
\times X = (Y \times S^1) \times X$ then $\mathcal{N}^3 = Y \times S^1$.\\

{\underline{{\tmstrong{Proof of the Claim}}}}: Let $\mathcal{M}^3 = Y \times
S^1$. Let $\beta \in \pi_1 (X)$ be a fixed non-trivial element and $\alpha \in
\pi_1 (\mathcal{N}^3)$ be an arbitrary element. We recall ({\it cf.} [T]) that the
centralizer of each nontrivial element of $\pi_1 (X)$ is an infinite cyclic
subgroup.

Now the homeomorphism $g : \mathcal{N}^3 \times X \longrightarrow
\mathcal{M}^3 \times X$ induces an isomorphism $g_{\ast} : \pi_1
(\mathcal{N}^3 \times X) \longrightarrow \pi_1 (\mathcal{M}^3 \times X)$ with
$g_{\ast} (\alpha, 1) = (\alpha', \alpha'')$ and $g_{\ast} (1, \beta) =
(\beta', \beta'')$.

Denote by $Z (\beta'')$ the centralizer of $\beta''$ and let $\gamma$ be its
generator, i.e., $\langle \gamma \rangle \cong Z (\beta'')$.

Since $(\alpha, 1) (1, \beta) = (1, \beta) (\alpha, 1)$, then $(\alpha',
\alpha'') (\beta', \beta'') = (\beta', \beta'') (\alpha', \alpha'')$. In
particular, $\alpha' \beta' = \beta' \alpha'$ and $\alpha'' \beta'' = \beta''
\alpha''$, and hence $\alpha'' \in Z (\beta'')$. Therefore there exist
integers $m, n$ such that $\alpha'' = \gamma^n$ and $\beta'' = \gamma^m$.
Since $\pi_1 (\mathcal{N}^3)$ is finitely generated (let us say by $\alpha_1,
\ldots, \alpha_k$ with $g_{\ast} (\alpha_i, 1) = (\alpha_i', \alpha_i'')$)
then $\alpha_i'' = \gamma^{n_i}$ for some $i = 1, \ldots, k$.

We claim that $\gamma^{n_i} = 1$ for all $i = 1, \ldots, k$; {\it i.e.\/}, $n_i =
0$, for $i = 1, 2, \ldots, k$.

To see this we argue by contradiction. Suppose the contrary. Then all
$\gamma^{n_i}$ generate an infinite cyclic subgroup of $Z (\beta'')$, namely a
subgroup generated by $\gamma^r$ where $r = \gcd (n_1, \ldots, n_k)$. In this
case we have an isomorphism
\[ \pi_1 (\mathcal{N}^3) \cong g_{\ast} (\pi_1 (\mathcal{N}^3), 1) \cong G
   \times Z \]
where $G \subseteq \pi_1 (\mathcal{M}^3)$ and $Z \subseteq \pi_1 (X)$. This
forces $\mathcal{N}^3
= Y \times S^1$ ({\it cf.} [He], p. 114). Suppose then that all $n_i = 0$. Then
\[ g_{\ast} (\pi_1 (\mathcal{N}^3), 1) \subseteq \pi_1 (\mathcal{M}^3) \times
   1 \]
If the above inclusion is proper then
\[ \pi_1 (X) \cong \pi_1 (\mathcal{M}^3) / p_1 g_{\ast} (\pi_1
   (\mathcal{N}^3), 1) \times \pi_1 (X) \]
where $p_1 : \pi_1 (\mathcal{M}^3) \times \pi_1 (X) \longrightarrow \pi_1
(\mathcal{M}^3)$ is the projection. This however is impossible. Consequently
$g_{\ast} (\pi_1 (\mathcal{N}^3), 1) = \pi_1
(\mathcal{M}^3)$, and hence $\pi_1 (\mathcal{N}^3) \cong \pi_1
(\mathcal{M}^3)$ and $\mathcal{N}^3 =\mathcal{M}^3 = Y \times S^1$.

This finishes the case when the genus of $S_1$ and $S_2$ is at least 2.

Finally, if one of $S_i, i = 1, 2$ (or both) is a torus, then the center of
one or both of $\pi_1 (\mathcal{M}^3)$ and \ $\pi_1 (\mathcal{N}^3)$ contains
at least two copies of $\mathbb{Z}$. This implies that one or both of
$\mathcal{M}^3$, $\mathcal{N}^3$ must be the torus $S^1 \times S^1 \times S^1$
and we are done.

Let us consider now the case of $\mathcal{M}^3$ and $\mathcal{N}^2$. Suppose
$\mathcal{M}^3 \times \mathcal{N}^2$ is decomposable i.e.
\[ \mathcal{M}^3 \times \mathcal{N}^2 = S_1 \times S_2 \times S^1 \]
where $S_1$, $S_2$ are surfaces.

Our considerations are divided into two cases:

(a) $\mathcal{N}^2 = S^2$.\\

(b) $\mathcal{N}^2$ has genus $> 1$.\\

\underline{Case (a)}: It follows that one of the $S_i, i = 1, 2$ must be
$S^2$, say $S_1 = S^2$. This gives $\mathcal{M}^3 \times S^2 = (S_2 \times
S^1) \times S^2$ and then $\mathcal{M}^3 = S_2 \times S^1$ because $\pi_1 (\mathcal{M}^3) \cong \pi_1 (S_2)\oplus \mathbb{Z}$.

\underline{Case (b)}: In analogy with the case of $\mathcal{M}^3,
\mathcal{N}^3$ it follows that, say $S_1 = \mathcal{N}^2$ and hence
\[ \mathcal{M}^3 \times \mathcal{N}^2 = (S_1 \times S^1) \times \mathcal{N}^2
\]
which implies $\mathcal{M}^3 = S_1 \times S^1$.

Finally we are left with the case $\mathcal{M}^3, \mathcal{N}^1$, so that
$\mathcal{N}^1 = S^1$. Then $\mathcal{M}^3 \times S^1 = S_1 \times S_2$. It
follows that say $S_1 = S^1 \times S^1$ and we can assume
$S_2 \neq S^1 \times S^1$. Indeed, if $S_1=S_2=S^1 \times S^1$ then
$\mathcal{M}^3=S^1 \times S^1 \times S^1$. Now
\[ \mathcal{M}^3 \times S^1 = (S^1 \times S_2) \times S^1 \]\\
Using the torus trick once again, we can arrange
$\mathcal{M}^3 \times \mathbb{R} = (S^1 \times S_2) \times \mathbb{R}$ and
hence $\pi_1(\mathcal{M}^3)~\cong~\pi_1(S_2)\;\oplus\;\mathbb{Z}$. This again
implies $\mathcal{M}^3 = S^1 \times S_2$ contradicting the indecomposability
of $\mathcal{M}^3$.

\

\

{\tmstrong{\underline{Proof of Theorem B}}}: In our proof we use a
$4$-manifold first constructed by S. Weinberger in [We] (see also Theorem 2.1
in [KS1]).  For completeness of this paper, we include a brief sketch of the
construction with somewhat different reasoning.  Let $\Sigma^3$ be a Seifert
homology 3-sphere with a natural free involution ({\it i.e.\/}, free
$\mathbb{Z}_2$-action) and Rochlin invariant ${\mu}
(\Sigma^3) = 1$. For example, we can let $\Sigma^3 = \{\Sigma (5, 7, 11),
\Sigma (3, 5,
13), \Sigma (3, 7, 17), \Sigma (5, 7, 27) \tmop{etc.} \ldots\}$ ({\it cf.}
[NR]).

Now $\Sigma^3 / \mathbb{Z}_2$ is a $\mathbb{Z}$-homology
$\mathbb{R}\mathbb{P}^3$, and there is a $\mathbb{Z} [\mathbb{Z}_2]$-homology
equivalence (see [KL], p. 35)
\[ f : \Sigma^3 / \mathbb{Z}_2 \longrightarrow \mathbb{R} \mathbb{P}^3~. \]
Let $I$ denote the interval [0,1], and
consider the map $h = f \times \tmop{id}_{I} : \Sigma^3 / \mathbb{Z}_2
\times I \longrightarrow \mathbb{R} \mathbb{P}^3 \times I$.
This map $h$ is a
$\mathbb{Z} [\mathbb{Z}_2] = \mathbb{Z} [\pi]$-homology
equivalence. If $\Gamma_0 (F)$, for $F = \tmop{id} : \mathbb{Z} [\pi]
\rightarrow \mathbb{Z} [\pi]$, is the Cappell-Shaneson homological surgery
group ({\it cf.} [CS]) then obviously $\lambda (h) = 0$ in $\Gamma_0 (F)$, here
$\lambda (h)$ is the surgery obstruction associated with $h$. But $L_0^{s, h}
(\pi) \cong \Gamma_0 (F)$ ({\it cf.} [CS], p. 289) and hence $\lambda (h) = 0$ in
$L_0^{s, h} (\pi)$.

Consequently
\[ h : \Sigma^3 / \mathbb{Z}_2 \times I \longrightarrow \mathbb{R}\mathbb{P}^3
   \times I \]
is normally bordant to a homotopy equivalence (rel boundary). By identifying the corresponding boundaries (using the identity mapping) we
obtain a homotopy equivalence $\mathcal{M} \longrightarrow
\mathbb{R}\mathbb{P}^3 \times S^1$ which we shall also call $h$.\\

\underline{Claim 1}: If $\widetilde{\mathcal{M}^{\nosymbol}}$ is the infinite
cyclic covering then there is no closed 3-manifold $\mathcal{N}^3$ with
$\widetilde{\mathcal{M}^{\nosymbol}} = \mathcal{N}^3 \times \mathbb{R}$.

\underline{Proof of Claim 1}: Suppose $\widetilde{\mathcal{M}^{\nosymbol}} =
\mathcal{N}^3 \times \mathbb{R}$. Then there is not difficult to see that,
there is a copy of $\Sigma^3 / \mathbb{Z}_2$ far away in the
$\mathbb{R}$-direction, which is disjoint with say $\mathcal{N}^3 \times
\{0\}$ in $\mathcal{N}^3 \times \mathbb{R}$.

The region between $\mathcal{N}^3 \times \{0\}$ and embedded $\Sigma^3 /
\mathbb{Z}_2$ is a homological $\mathbb{Z} [\mathbb{Z}_2]$ $h$-cobordism
$(\mathcal{W} ; \mathcal{N}^3 ; \Sigma^3 / \mathbb{Z}_2)$. Since
$\mathcal{N}^3$ and $\Sigma^3 / \mathbb{Z}_2$ are parallelizable there is a
prefered spin structure on $\mathcal{N}^3$ and $\Sigma^3 / \mathbb{Z}_2$. One
can ask about the possibility of extending this structure to the entire
$\mathcal{W}$. Whether one can do this or not is determined by the obstruction
in $H^4 (\mathcal{W} ; \partial \mathcal{W} ; \mathbb{Z}_2) \simeq
\mathbb{Z}_2$ ({\it cf.} [KS1], p. 448).

There is an analogous obstruction for the existence of a spin structure on the
2-fold cover $(\widetilde{\mathcal{W}^{\nosymbol} ;} \partial
\widetilde{\mathcal{W}^{\nosymbol}})$. By the naturality, the obstruction in
$H^4 (\widetilde{\mathcal{W}^{\nosymbol} ;} \partial
\widetilde{\mathcal{W}^{\nosymbol}} ; \mathbb{Z}_2) \simeq \mathbb{Z}_2$ is
the image of the obstruction in $H^4 (\mathcal{W} ; \partial \mathcal{W} ;
\mathbb{Z}_2)$ under the homomorphism
\[ H^4 (\mathcal{W} ; \partial \mathcal{W} ; \mathbb{Z}_2) \longrightarrow
   H^4 (\widetilde{\mathcal{W}^{\nosymbol} ;} \partial
   \widetilde{\mathcal{W}^{\nosymbol}} ; \mathbb{Z}_2) \]
and hence the obstruction in $H^4 (\widetilde{\mathcal{W}^{\nosymbol} ;}
\partial \widetilde{\mathcal{W}^{\nosymbol}} ; \mathbb{Z}_2)$ is trivial. In
other words the manifold $(\widetilde{\mathcal{W}^{\nosymbol} ;} \partial
\widetilde{\mathcal{W}^{\nosymbol}})$ is quasi-spin in the terminology of
[KS1], p. 449.

In our case $(\widetilde{\mathcal{W}^{\nosymbol} ;} \partial
\widetilde{\mathcal{W}^{\nosymbol}})$ is a $\mathbb{Z}$-homological
$h$-cobordism between $\Sigma^3$ and a homotopy 3-sphere
$\widetilde{\mathcal{N}^3}$. The Rochlin ${\mu}$-invariant is invariant
with respect to topological quasi-spin $\mathbb{Z}$-homological-$h$-cobordism ([KS1]) and
hence ${\mu} (\widetilde{\mathcal{N}^3}) ={\mu} (\Sigma^3)$. 
However, this is a contradiction since ${\mu} (\Sigma^3) = 1$ and ${\mu}
(\widetilde{\mathcal{N}^3}) = 0$ by the Casson's results ({\it cf.} [AM]), and hence
Claim 1 has been established.

\

Now let $h : \mathcal{M} \rightarrow \mathbb{R}\mathbb{P}^3 \times S^1$ be the
constructed homotopy equivalence.\\

\underline{Claim 2}: The homotopy equivalence $h : \mathcal{M} \rightarrow
\mathbb{R}\mathbb{P}^3 \times S^1$ is normally bordant to the identity.

\underline{Proof of Claim 2}: Consider the Wall-Sullivan exact surgery
sequence ({\it cf.} [Wa]), which extends to dimension 4 by the results of [FQ]:
\[ \ldots \longrightarrow L^s_1  (\mathbb{Z} \times \mathbb{Z}_2)
   \overset{\gamma}{\longrightarrow} S^{\tmop{TOP}}  (S^1 \times \mathbb{R}
   \mathbb{P}^3) \overset{\eta}{\longrightarrow} [S^1 \times \mathbb{R}\mathbb{P}^3 ; G /
   \tmop{TOP}] \overset{\lambda}{\longrightarrow} L^s_0  (\mathbb{Z} \times
   \mathbb{Z}_2) \]
Here $[S^1 \times \mathbb{R}\mathbb{P}^3 ; G / \tmop{TOP}] \cong H^2  (S^1 \times
\mathbb{R}\mathbb{P}^3 ; \mathbb{Z}_2) \oplus H^4  (S^1 \times \mathbb{R}\mathbb{P}^3 ;
\mathbb{Z}) \cong \mathbb{Z}_2 \oplus \mathbb{Z}_2 \oplus \mathbb{Z}$.

By [Sh] and the triviality of $\mathrm{Wh}(\mathbb{Z}_2)$ we have
\[ L^s_1  (\mathbb{Z} \times \mathbb{Z}_2) \cong L^s_1 (\mathbb{Z}_2)
   \oplus L^s_0 (\mathbb{Z}_2) \]
and
\[ L^s_0  (\mathbb{Z} \times \mathbb{Z}_2) \cong L^s_0 (\mathbb{Z}_2)
   \oplus L^s_3 (\mathbb{Z}_2) \cong L^s_0 (0) \oplus \widetilde{L_0^s}
   (\mathbb{Z}_2) \oplus \mathbb{Z}_2 \]
in which $L^s_0 (0)$ and $\widetilde{L_0^s} (\mathbb{Z}_2)$ are both
isomorphic to $\mathbb{Z}$.

Let us briefly analyze the subgroup $\tmop{Ooze} (\mathbb{Z} \times
\mathbb{Z}_2) \subset L^s_0  (\mathbb{Z} \times \mathbb{Z}_2)$. We recall
that the $\tmop{Ooze} (-)$ subgroup is represented by surgery obstructions
determined by closed manifolds ({\it cf.} [HMTW]).

It turns out that $\tmop{Ooze} (\mathbb{Z} \times \mathbb{Z}_2) \cong L^s_0
(0) \oplus L^s_3 (\mathbb{Z}_2) \cong \mathbb{Z} \oplus \mathbb{Z}_2$.To
see this, just note that $\mathbb{Z} \cong L^s_0 (0)$ is represented by the
difference of signatures and the existence of the $E_8$ manifold ({\it cf.} [FQ])
implies $L^s_0 (0) \subset \tmop{Ooze} (\mathbb{Z} \times \mathbb{Z}_2)$.

The copy of $\mathbb{Z}_2 \cong L^s_3 (\mathbb{Z}_2)$ is determined by the
codimension one Arf invariant ({\it cf.} [HMTW], Theorem A) in $[\mathbb{R}\mathbb{P}^3;
\mathbb{Z}_2] \overset{\lambda}{\longrightarrow} L^s_3 (\mathbb{Z}_2)$.

In fact $\mathbb{Z}_2 \cong [\mathbb{R}\mathbb{P}^3 ; G / \tmop{TOP}] \cong H^2
(\mathbb{R}\mathbb{P}^3 ; \mathbb{Z}_2)$ corresponds to the copy of $\mathbb{Z}_2
\subset H^2  (S^1 \times \mathbb{R}\mathbb{P}^3 ; \mathbb{Z}_2)$ given by $H^0 (S^1
; \mathbb{Z}_2) \otimes H^2  (\mathbb{R}\mathbb{P}^3 ;
\mathbb{Z}_2)$.

The remaining copy of $\mathbb{Z}_2 \subset H^2  (S^1 \times \mathbb{R}\mathbb{P}^3
; \mathbb{Z}_2)$ corresponds to $H^1 (S^1 ; \mathbb{Z}_2) \otimes H^1
(\mathbb{R}\mathbb{P}^3 ; \mathbb{Z}_2)$.

It turns out that this remaining copy of $\mathbb{Z}_2$ is represented by a
self-homotopy equivalence
\[ s : S^1 \times \mathbb{R}\mathbb{P}^3 \longrightarrow S^1 \times \mathbb{R}\mathbb{P}^3
\]
given by the pinching construction. To be more specific the self-homotopy
equivalence $s$ is given by
\[ s : S^1 \times \mathbb{R}\mathbb{P}^3 \overset{\vee}{\longrightarrow} S^1 \times
   \mathbb{R}\mathbb{P}^3 \vee S^4 \overset{\tmop{id} \vee v}{\longrightarrow} S^1
   \times \mathbb{R}\mathbb{P}^3 \]
where $v : S^4 \rightarrow S^1 \times \mathbb{R}\mathbb{P}^3$ is the nontrivial homotopy
class of $\pi_4 (S^3) \cong \mathbb{Z}_2$.

The homotopy theoretic argument for this is given in [KS2], Theorem 2.1, Case II,
p. 531.

The argument in [KS2] is given for a self-homotopy equivalence (rel boundary)
of $\mathbb{R}\mathbb{P}^3 \times I$, but it works in precisely the same way for $S^1
\times \mathbb{R}\mathbb{P}^3$.

A consequence of the above considerations is that the homotopy (simple)
equivalence
\[ h : \mathcal{M} \longrightarrow S^1 \times \mathbb{R}\mathbb{P}^3 \]
being an element of $S^{\tmop{TOP}}  (S^1 \times \mathbb{R}\mathbb{P}^3)$ must be
normally bordant to the identity.

(We recall that $\tmop{Wh} (\mathbb{Z} \times \mathbb{Z}_2) \cong \tmop{Wh}
(\mathbb{Z}_2) \oplus \widetilde{K_0} (\mathbb{Z} [\mathbb{Z}_2]) \cong
0$; see [BHS] and [Ha]).

Let $(\mathcal{W} ; \mathcal{M}, S^1 \times \mathbb{R}\mathbb{P}^3)$ be the
corresponding normal bordism:

\begin{figure}
\centering
\includegraphics[width=0.9\textwidth]{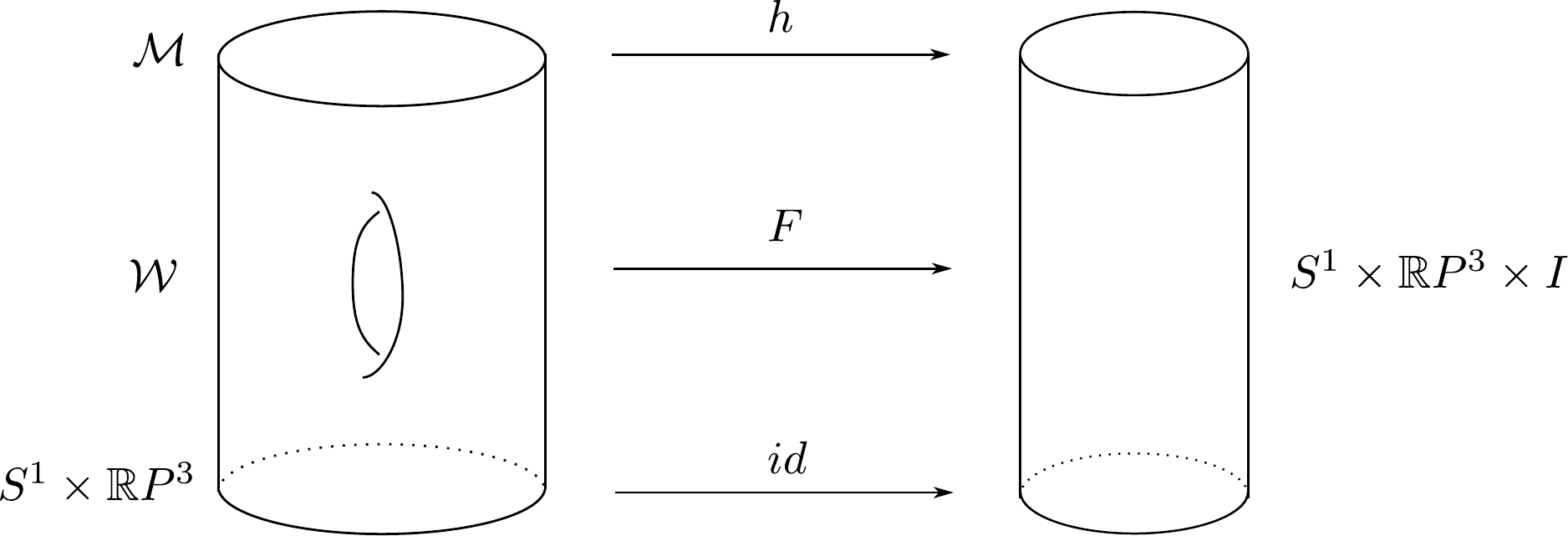}
\caption{The normal cobordism $(\mathcal{W} ; \mathcal{M}, S^1 \times \mathbb{R}\mathbb{P}^3)$.}
\label{fig:fig1}
\end{figure}

Multiplying the above normal bordism by the $\tmop{id}_{S^i} : S^i \rightarrow
S^i~(i = 2, 3, 4)$ we get the surgery obstruction ({\it cf.} [Mo]):
\[ \lambda (F \times \tmop{id}_{S^i}) = \lambda (F) \cdot \sigma^{\ast} (S^i)
   = 0 \].
In particular $S^1 \times \mathbb{R}\mathbb{P}^3 \times S^i$ is s-cobordant to
$\mathcal{M} \times S^i$ and hence $\mathcal{M} \times S^i = S^1 \times
\mathbb{R}\mathbb{P}^3 \times S^i$.

This finishes the proof of Theorem B, once we know that the manifold
$\mathcal{M}$ is indecomposable. This however follows from the Claim 1. We are then
left with the construction of infinitely many corresponding examples.

To do this we follow [KS1]. Consider the extension of the
Wall-Sullivan exact sequence to dimensional 3 ({\it cf.} [JK]):
\[ \ldots \longrightarrow L^s_0 (\mathbb{Z}_2) \longrightarrow S^H
   (\mathbb{R}\mathbb{P}^3) \longrightarrow [\mathbb{R}\mathbb{P}^3 ; G / \tmop{TOP}]
   \overset{}{\longrightarrow} L^s_3 (\mathbb{Z}_2) \]
Now $L^s_0 (\mathbb{Z}_2) \cong L^s_0 (0) \oplus \widetilde{L_0^s}
(\mathbb{Z}_2) \cong \mathbb{Z} \oplus \mathbb{Z}$ and $\widetilde{L_0^s}
(\mathbb{Z}_2)$ acts freely on $S^H  (\mathbb{R}\mathbb{P}^3)$. This implies

\

(a)~the existence of infinitely many homology 3-spheres $(\Sigma_i^3, t_i), i =
1, 2, \ldots$ with ${\mu} (\Sigma_i^3) = 1$ and a free involution $t_i :
\Sigma_i \rightarrow \Sigma_i$

(b)~the $\rho$-invariant associated with these actions are different, i.e.
$\rho (\Sigma_i^3, t_i) - \rho (\Sigma_j^3, t_j) \neq 0$ for $i \neq j$.

\

The crucial fact needed here is the congruence
\[ {\mu} (\Sigma^3_i) \equiv \rho (\Sigma_i^3, t_i) \tmop{mod} 16 \]
({\it cf.} [NR]).

Given the above we start with $\Sigma^3_i / \mathbb{Z}_2 \times I$ and
convert it by topological surgery to a homotopy equivalence
\[ (\overline{\mathcal{W}_i}, \partial) \longrightarrow (\mathbb{R}\mathbb{P}^3
   \times I, \partial) \]
which is a $\mathbb{Z} [\mathbb{Z}_2]$-homology equivalence of boundaries.

Next form a two ended open manifold $\widetilde{\mathcal{W}_i}$ by taking
infinitely many copies of $\overline{\mathcal{W}_i}$, one on the top of the another.

\begin{figure}
	\centering
	\includegraphics[width=0.9\textwidth]{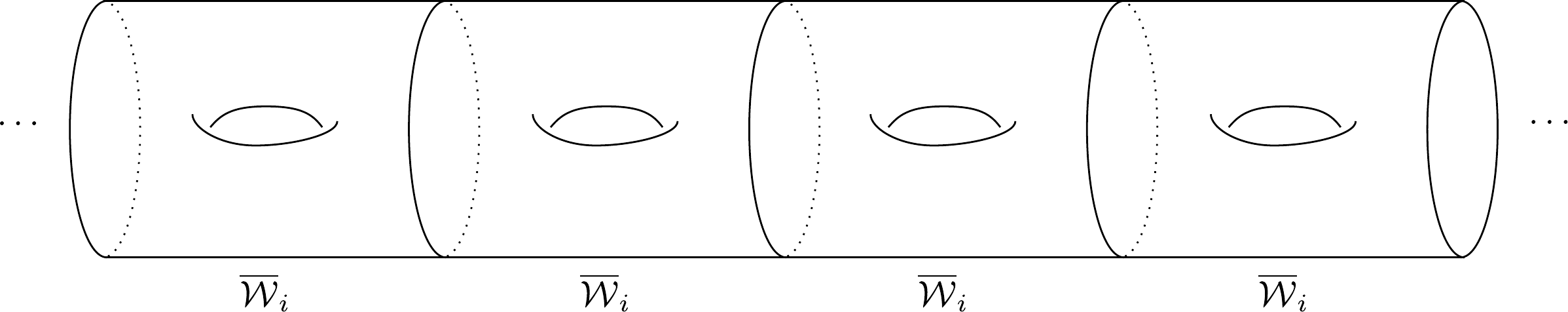}
	\caption{Open manifold $\widetilde{\mathcal{W}_i}$.}
	\label{fig:fig2}
\end{figure}

There is a natural free proper action of $\mathbb{Z}$ on
$\widetilde{\mathcal{W}_i}$ and we shall let
\[ \mathcal{W}_i \assign \widetilde{\mathcal{W}_i} / \mathbb{Z}~. \]
The manifold $\mathcal{W}_i~(i = 1, 2, \ldots)$ has the required properties,
more precisely, we have the following:

(1)~$\mathcal{W}_i \simeq S^1 \times \mathbb{R}\mathbb{P}^3~(i = 1, 2, \ldots)$.

(2)~$\mathcal{W}_i$ is indecomposable.

(3)~$\mathcal{W}_i \neq \mathcal{W}_j$, $i \neq j$.

(4)$\mathcal{W}_i \times S^k = S^1 \times \mathbb{R}\mathbb{P}^3 \times S^k~(k = 2, 3,
4)$.

\

\

\underline{{\tmstrong{Proof of Theorem B$'$}}}: Suppose $\mathcal{M}^4 \times
S^4$ is decomposable. Then there are two cases:

(a)~$\mathcal{M}^4 \times S^4 = S^2 \times K^6$

(b)~$\mathcal{M}^4 \times S^4 = \Sigma^3 \times K^5$

for some closed manifolds $K^6$ and $K^5$, and a homotopy 3-sphere $\Sigma^3$.
We know that $\Sigma^3 = S^3$ by Perelmen's proof of the Poincar\'e
Conjecture ({\it cf.} [DL]), but we do not need this result.

The case (b) cannot occur because $\chi (\mathcal{M}^4 \times S^4) \neq 0$ and $\chi
(\Sigma^3 \times K^5) = 0$.

Consider now case (a). Note that one can assume $H^2 (\mathcal{M}^4) \neq 0,
\mathbb{Z}$. Suppose then that $H^2 (\mathcal{M}^4) \cong \mathbb{Z} \oplus
\mathbb{Z}$. Then
\[ H^2  (\mathcal{M}^4 \times S^4) \cong H^2 (\mathcal{M}^{\nosymbol}) \cong
   H^2 (S^2) \oplus H^2 (K^6) \]
It follows that the matrix for the standard intersection form on $H^2
(\mathcal{M})$ is given by $\left( \begin{array}{cc}
  0 & 1\\
  1 & 0
\end{array} \right)$. This however implies $\mathcal{M} = S^2 \times S^2$
which is a contradiction.

Next, suppose $H^2 (\mathcal{M}^4) \cong \underset{r}{\oplus} \mathbb{Z}$,
where $r \geqslant 3$. In this case the matrix for the intersection form \ on $H^2
(\mathcal{M}^4)$ is given by
$$\left( \begin{array}{ccccc}
  0 & 1 & 1 & \ldots & 1\\
  1 & 0 & 0 & \cdots & 0\\
  1 & 0 & 0 &  & \\
  \vdots & \vdots &  & \ddots & \\
  1 & 0 &  &  & 0
\end{array} \right)~.$$
This again is a contradiction because the determinant of
the displayed matrix is zero).

The case of $\mathcal{M}^4 \times S^3$ is easier. The only possible
decomposition of $\mathcal{M}^4 \times S^3$ is given by $\mathcal{M}^4 \times
S^3 = \Sigma^3 \times S^2 \times S^2$ for a homotopy 3-sphere $\Sigma^3$.

Finally, $\mathcal{M}^4 \times S^2$ can only be decomposed as $\mathcal{M}^4
\times S^2 = S^2 \times S^2 \times S^2$.

\
\\

{\noindent}\tmtextbf{Acknowledgments . } We are grateful to Professor Witold
Rosicki for numerous discussions and suggestions which have improved this
paper in several ways.{\medskip}
\\
\\
\\

\begin{flushleft}
{\Large \textbf{References}}
\end{flushleft}

\begin{flushleft}
[AM]\quad S. Akbulut and J. McCarthy, Casson's invariant for oriented homology
3-spheres: an exposition, Princeton University Press, Princeton, 1990.

\

[BHS]\quad H. Bass, A. Heller and R. Swan, The Whitehead group of polynomial extension,
Publ. Inst. Hautes \'Etudes Sci. 22 (1964), 35--63.

\

[Bi]\quad R. H. Bing, The cartesion product of a certain non-manifold and a line is
$E^4$, Ann. Math. 70 (1959), 399--412.

\

[Bo]\quad K. Borsuk, On decomposition of manifolds into products of curves and
surfaces, Fund. Math. 33 (1945), 273--298.

\

[CJ]\quad A. Casson and D. Jungreis, Convergence groups and Seifert fibered
3-manifolds, Invent. Math. 118 (1994), 441--456.

\

[CR]\quad P. Conner and F. Raymond, Derived actions in Proc. Second Conf. Compact
Transf. Groups, Amherst, 1971, in: Lecture Notes in Math. vol. 299 Springer,
Berlin/New York, 1972, pp. 237--310.

\

[CS]\quad S. Cappell and J. Shaneson, The codimension two placement problem and
homology equivalent manifolds, Ann. Math. 99 (1974), 277--348.

\

[FJ]\quad F. T. Farrell and L. E. Jones, Topological rigidity for compact nonpositively
curved manifolds, Proc. Sympos. Pure Math. 54 Part 3 (1993), pp. 229--274.

\

[FQ]\quad M. Freedman and F. Quinn: topology of 4-manifolds, Princeton Math. Series
39, Princeton Univ. Press, Princeton, 1990.

\

[G]\quad D. Gabai, Convergence groups are Fuchsian groups, Ann. Math. 136(1992),
447--510.

\

[Ha]\quad D. R. Harmon, $NK_1$ of finite groups, Proc. Amer. Math. Soc. vol. 100, No.2
(1987) 229--232.

\

[He]\quad J. Hempel, 3-manifolds, Annals of Mathematics Studies, No.86, Princeton
Univ. Press, Princeton, 1976.

\

[HMTW]\quad I. Hambleton, R. J. Milgram, L. Taylor and B. Williams, Surgery with finite
fundamental group, Proc. London Math. Soc. (3) 56 (1988), 349--379.

\

[JK]\quad B. Jahren and S. Kwasik, Three-dimensional surgery theory, UNil-groups and
the Borel Conjecture, Topology 42 (2003), 1353--1369.

\

[JS]\quad W. H. Jaco and P. Shalen, Seifert fibered spaces in 3-manifolds, Memoirs of
the AMS, vol 21 (1979) No. 220.

\

[KL]\quad S.Kwasik and T.Lawson, Nonsmoothable $\mathbb{Z}_p$ actions on
contractible 4-manifolds, J. reine angew. Math. 437 (1993), 29--54.

\

[KR1]\quad S.Kwasik and W.Rosicki, Cartesian product stabilization of 3-manifolds,
Topology Appl. 157 (2010), 2342--2346.

\

[KR2]\quad S.Kwasik and W.Rosicki, On stability of 3-manifolds, Fund. Math. 
182 (2004), 169--180.

\

[KS1]\quad S. Kwasik and R. Schultz, Desuspension of group actions and the Ribbon
Theorem, Topology 27 (1988), 443--457.

\

[KS2]\quad S. Kwasik and R. Schultz, On $s$-cobordisms of metacyclic prism manifolds,
Invent. Math. 97 (1989), 523--552.

\

[M]\quad J. Morgan, A product formula for surgery obstructions, Mem. of the Amer.
Math. Soc. Vol.14, No. 201, AMS, 1978.

\

[NR]\quad W. Neumann and F. Raymond, Seifert manifolds, plumbing ${\mu}$-invariant
and orientation reversing maps, Lecture Notes in Math. 664, 163--169,
Springer-Verlag, 1978.

\

[R]\quad D. Rolfsen, Knots and Links, Math. Lecture Ser. Vol. 7, Publish or Perish,
Berkeley 1976.

\

[Sh]\quad J. Shaneson, Wall's surgery obstruction group for $\mathbb{Z} \times G$,
Ann. Math. 90 (1969), 296--334.

\

[T]\quad B. Truffault, Centralisateurs des elements dans les groups de Greendlinger.
C. R. Acad. Sci. Paris, Ser. A. 279 (1974). pp. 317--319.

\

[Wa]\quad C. T. C. Wall, Surgery on Compact Manifolds, London Math. Soc. monographs,
No. 1, Academic Press, London and New York, 1970.

\

[We]\quad S. Weinberger, On fibering four- and five-manifolds, Isr. J. Math.
59 (1987), 1-7.
\end{flushleft}

\

\

\

\

\

\

\

\end{document}